\documentclass[11pt]{article}
\usepackage{amssymb}
\usepackage{mathrsfs}
\usepackage{latexsym,amsmath,amssymb,amsfonts,epsfig,graphicx,cite,psfrag}
\usepackage{eepic,color,colordvi,amscd}
\usepackage{ebezier}
\usepackage{verbatim}
\usepackage{comment}

\newcommand{\qed}{\hfill $\Box $}
\newcommand{\pf}{\noindent {\bf Proof.} }

\newtheorem{theorem}{Theorem}[section]
\newtheorem{lemma}[theorem]{Lemma}
\newtheorem{coro}[theorem]{Corollary}
\newtheorem{proposition}[theorem]{Proposition}
\newtheorem{conjecture}[theorem]{Conjecture}

\newtheorem{problem}[theorem]{Problem}
\newtheorem{question}[theorem]{Question}

\begin{document}

\title{On rainbow matchings for hypergraphs}

\author{Hongliang Lu \footnote{luhongliang@mail.xjtu.edu.cn; partially supported by the National Natural
Science Foundation of China under grant No.11471257 and
Fundamental Research Funds for the Central Universities}
\\School of Mathematics and Statistics\\
Xi'an Jiaotong University\\
Xi'an, Shaanxi 710049, China\\
\smallskip\\
Xingxing Yu \footnote{yu@math.gatech.edu; partially supported by NSF grants DMS-1265564 and DMS-1600387}\\
School of Mathematics\\
Georgia Institute of Technology\\
Atlanta, GA 30332, USA\\}

\date{}
\maketitle

\begin{abstract}
%Let $n,k,t$ be positive integers. For $\mathcal{F}=\{\mathcal{F}_1,\ldots,\mathcal{F}_t\}$ with each ${\cal
%  F}_i$ a family of subsets of a given set $S$, a rainbow matching
%$M$ for $\mathcal{F}$ is a collection of  $t$ pairwise disjoint
%subsets of $S$ such that $|M\cap \mathcal{F}_i|=1$ for $i=1, \ldots,
%t$.
For any posotive integer $m$, let $[m]:=\{1,\ldots,m\}$. 
Let $n,k,t$ be positive integers. Aharoni and Howard conjectured that
if, for $i\in [t]$,
$\mathcal{F}_i\subset[n]^k:= \{(a_1,\ldots,a_k): a_j\in [n] \mbox{ for
} j\in [k]\}$ and $|\mathcal{F}_i|>(t-1)n^{k-1}$, then there exist
$M\subseteq [n]^k$ such that $|M|=t$ and $|M\cap \mathcal{F}_i|=1$ for
$i\in [t]$
%then $\{\mathcal{F}_1,\ldots,\mathcal{F}_t\}$ admits a rainbow matching. 
We show that this conjecture holds when $n\geq 3(k-1)(t-1)$.

Let $n, t, k_1\ge k_2\geq \ldots\geq k_t $ be positive integers.  Huang, Loh and Sudakov asked for the maximum  $\Pi_{i=1}^t |{\cal
   R}_i|$ over all ${\cal R}=\{{\cal R}_1, \ldots ,{\cal R}_t\}$  such that
each ${\cal R}_i$ is a collection of $k_i$-subsets of $[n]$ for which 
there does not exist a collection $M$ of subsets of $[n]$ such that $|M|=t$ and $|M\cap \mathcal{R}_i|=1$ for $i\in [t]$
%and ${\cal R}$ does not admit a rainbow matching.
We show that  for  sufficiently large $n$ with $\sum_{i=1}^t k_i\leq n(1-(4k\ln n/n)^{1/k}) $,
$\prod_{i=1}^t |\mathcal{R}_i|\leq {n-1\choose k_1-1}{n-1\choose
  k_2-1}\prod_{i=3}^{t}{n\choose k_i}$. This bound is tight. 
\end{abstract}

\section{Introduction}

For a positive integer $k$ and a set $S$, let $[k]:=\{1,\ldots, k\}$ and
${S\choose k}:=\{T\subseteq S: |T|=k\}$. A {\it hypergraph} $H$
consists of a vertex set $V(H)$ and an edge set $E(H)\subset 2^{V(H)}$.
% whose members are subsets of $V(H)$.
Thus, for any positive integer $n$, any subset of $2^{[n]}$ forms a hypergrpah with vertex set $[n]$.

Let $k$ be a positive integer. A hypergraph $H$ is {\it $k$-uniform} if $E(H)\subseteq {V(H)\choose k}$, and a $k$-uniform hypergraph is also
called a {\it $k$-graph}. A $k$-graph $H$ is \emph{$k$-partite }if there exists a
partition of $V(H)$ into sets $V_1, \cdots, V_k$ (called {\it partition  classes}) such that for any $f\in E(H)$, $|f\cap V_i|=1$
for $i\in [k]$.

Let $H$ be a hypergraph and $T\subseteq V(H)$. We write
$e(H):=|E(H)|=|H|$. (Note that we often identify $E(H)$ with $H$.) 
The {\it degree}
of $T$ in $H$, denoted by $d_H(T)$, is the number of edges of $H$
containing $T$. For any integer  $l\ge 0$, let $\delta_l(H):=
\min\{d_H(T): T\in {V(H)\choose l}\}$ denote the minimum
{\it $l$-degree} of $H$. Hence,  $\delta_0(H)$ is the number of edges in $H$.  Note that $\delta_1(H)$ is often called the minimum {\it
vertex} degree of $H$.  If $H$ is a $k$-graph then $\delta_{k-1}(H)$ is also known as the minimum
{\it codegree} of $H$.

Let $H$ be a $k$-partite $k$-graph, with partition classes $V_1,\ldots, V_k$.  We say that $H$ is {\it balanced} if $|V_i|=|V_j|$ for all $i,j\in [k]$.
A set $T\subseteq V(H)$ is said to be {\it
legal} if $|T\cap V_i|\leq 1$ for all $i\in [k]$.  Thus, if $T$ is
not legal in $H$ then $d_H(T)=0$. So for integer $l$ with $0\le l\le k-1$, let $\delta_l(H):= \min\{d_H(T):
T\in {V(H)\choose l} \mbox{ and $T$ is legal}\}$.

A \emph{matching} in a hypergraph $H$ is a set of pairwise disjoint
edges in $H$, and we use $\nu(H)$ to denote the maximum size of a
matching in $H$. A classical problem in extremal set
theory is to determine $max |H|$ with $\nu(H)$ fixed. Erd\H{o}s
\cite{Erdos65} in 1965 made the following conjecture:
For positive integers $k,n,t$, every $k$-graph $H$ on $n$ vertices with $\nu(H) <
t$ satisfies $e(H)\leq \max \left\{{kt-1\choose k}, {n\choose k}-{n-t+1\choose k}\right\}.$
This bound is tight because of the complete $k$-graph on $kt-1$ vertices and the $k$-graph on $n$ vertices in which every 
edge intersects a fixed set of $t$ vertices. 

There has been attempts to extend the above conjecture of  Erd\H{o}s to a family
of hypergraphs. Let $\mathcal{F} = \{{\cal F}_1,\ldots, {\cal F}_t\}$ be a family
of hypergraphs. A set of pairwise disjoint edges, one from
each ${\cal F}_i$, is called a \emph{rainbow matching} for $\mathcal{F}$. (In this case, we also say that ${\cal F}$ or
$\{{\cal F}_1,\ldots, {\cal F}_t\}$ {\it admits} a rainbow matching.)
 Huang, Loh and Sudakov \cite{HLS} and,   independently,  Aharoni and
 Howard \cite{AH} made the following conjecture:
%\begin{conjecture}\label{HLSAH}
Let $t$ be a positive integer and ${\cal F}=\{\mathcal{F}_1,\ldots,\mathcal{F}_t\}$ such that, for $i\in [t]$,
${\cal F}_i\subseteq {[n]\choose k}$ and
$|\mathcal{F}_i|> \left\{{kt-1\choose k}, {n\choose k}-{n-t+1\choose k}\right\};$
 then  ${\cal F}$ admits a rainbow matching.
%\end{conjecture}
Huang, Loh and Sudakov \cite{HLS} showed that this conjecture  holds for $n> 3k^2t$.
Aharoni and Howard  \cite{AH} also proposed the following $k$-partite
version.

\begin{conjecture}\label{AH-conj}
 If ${\cal F}=\{\mathcal{F}_1,\ldots,\mathcal{F}_t\}$ such that, for $i\in [t]$,  ${\cal
   F}_i$ is a $k$-partite $k$-graph in which each partition class has size $n$, 
%\subset [n]^k:=\{(a_1,\ldots,a_k): a_j\in [n] \mbox{ for } j\in [k]\}$ and $|{\cal F}_i|>(t-1)n^{k-1}$, then
then ${\cal F}$ admits a rainbow matching.
\end{conjecture}

 Aharoni and Howard  \cite{AH} proved Conjecture~\ref{AH-conj} for $t=2$ or $k\leq 3$. Our first result
 implies that  Conjecture~\ref{AH-conj} holds when $n\geq 3(k-1)(t-1)$.

\begin{theorem}\label{TH-AHC}
Let $k,r,n, t$ be positive integers such that $2\leq r\leq k$ and $n\geq
3(k-1)(t-1)$, and let $U_1, \ldots, U_k$ be pairwise disjoint sets with $|U_i|=n$ for $i\in [k]$.
 For each $i\in [t]$, let $i_1,\ldots, i_r\in [k]$ be pairwise distinct such that
$\mathcal {F}_i\subseteq U_{i_1}\times \cdots\times U_{i_r}$ and $|\mathcal {F}_i|>(t-1)n^{r-1}$.
Then  $\{{\cal F}_1,\ldots,{\cal F}_t\}$ admits a rainbow matching.
\end{theorem}

The famous Erd\H{o}s-Ko-Rado theorem states that if $k \leq n/2$ and 
$H\subset {[n]\choose k}$ has more than
${n-1\choose k-1}$ edges, then $\nu(H) > 1$.
%This has been extended in more than one way to pairs of hypergraphs.
Pyber \cite{PY} gave a  product-type generalization of the
Erd\"os-Ko-Rado theorem, which was improved by Matsumoto and Tokushige \cite{MN} to the following:
%\begin{theorem}[Matsumoto and Tokushige, \cite{MN}]\label{MTT}
Let $k_1,k_2,n$ be positive integers such that $n\geq \max\{2k_1,2k_2\}$,
and let $H_i\subset {[n]\choose k_i}$ for $i\in [2]$ such that  $e(H_1)e(H_2)>
{n-1\choose k_1-1}{n-1\choose k_2-1}$; then $\{H_1,H_2\}$ admits a rainbow
matching.
%\end{theorem}
Huang, Loh and Sudakov  \cite{HLS} asked the following more general question.

\begin{problem}\label{HLSP}
For positive integers $n,t,k_1, \ldots, k_t$, what is the maximum $\Pi_{i=1}^t |{\cal R}_i|$ among families
${\cal R}=\{{\cal R}_1, \ldots, {\cal R}_t\}$ such that ${\cal
  R}_i\subseteq {[n]\choose k_i}$ for $i\in [t]$ and ${\cal R}$ admits no rainbow matching.
\end{problem}

Our second result in this paper provides an answer to Problem
\ref{HLSP} when $n$ is large.

\begin{theorem}\label{TH-HLS}
  Let $n,t, k_1,\ldots, k_t$ be positive integers such that $n$ is sufficiently large
  and $\sum_{i=1}^t k_i\leq n(1-(\frac{8k\ln n}{n})^{1/k})$. Suppose $k_1,k_2\ge k_i$ for $i=3, \ldots, t$. 
  Let $\mathcal{F}_i\subset {[n]\choose k_{i}}$ for $i\in [t]$, such that
%Let $\pi$ be a permutation of $[t]$ and let $\mathcal{F}_i\subseteq {[n]\choose k_{\pi(i)}}$ for $i\in [t]$.
\[
 |\mathcal{F}_1|\ |\mathcal{F}_2| > {n-1\choose k_1-1}{n-1\choose k_2-1}\quad \mbox{when $t=2$,}
\]
and
\[
\prod_{i=1}^t |\mathcal{F}_i|> {n-1\choose k_1-1}{n-1\choose k_2-1}\prod_{i=3}^t{n\choose k_i}\quad \mbox{when $t\geq 3$}.
\]
Then  $\{{\cal F}_1,\ldots, {\cal F}_t\}$ admits a rainbow matching.
\end{theorem}

We remark that the bound in Theorem~\ref{TH-HLS} is tight.  
%Let $k_1\geq k_2\geq  \ldots\geq k_t\ge 2$ be integers. 
Let $\mathcal{F}_i=\{e \ :\ 1\in e \in  {[n]\choose k_i}\}$ for $i\in
[2]$ and let $\mathcal{F}_i={[n]\choose k_i}$ for $i\in
[t]-\{1,2\}$. Then
\[
\prod_{i=1}^t |\mathcal{F}_i|={n-1\choose k_1-1}{n-1\choose k_2-1}\prod_{i=3}^t {n\choose k_i}.
\]
Clearly,  $\{\mathcal{F}_1,\mathcal{F}_2\}$ does not admit  any rainbow
matching. Hence,  $\{\mathcal{F}_1,\ldots,\mathcal{F}_t\}$ does not
admit any  rainbow matching.

Our third result is a natural extension of Theorem 3.3 in \cite{HLS} by Huang, Loh and Sudakov.

\begin{theorem}\label{HLS-EX12}
Let $n,t, k_1,\ldots k_t$ be positive integers such that $n> 3k^2 t$,
and let $k=\max\{k_i: i\in [t]\}$.
For $i\in [t]$, let $\mathcal{F}_i\subset {[n]\choose k_i}$ such that $|\mathcal{F}_i|>{n\choose k_i}-{n-t+1\choose k_i}$.
Then $\{\mathcal{F}_1, \ldots, \mathcal{F}_t\}$ admits a rainbow matching.
\end{theorem}

In view of Theorem~\ref{HLS-EX12}, we ask the following 
% more specific version of Problem \ref{HLSP}.
\begin{question}
Let $k,n,t$ be positive integers and let $\varepsilon$ be a constant such that $0<\varepsilon<1$ and $n\ge kt/(1-\varepsilon)$, and let $\mathcal{R}_i\subset
{[n]\choose k}$ for $i\in [t]$ such that $|\mathcal{R}_1|\leq |{\cal
  R}_2|\le \ldots\leq |\mathcal{R}_t|$. Is it true that if for all $r\in [t]$,
\[
\prod_{i=1}^r |\mathcal{R}_i|>\left({n\choose k}-{n-r+1\choose k}\right)^r,
\]
 then $\{\mathcal{R}_1, \ldots, \mathcal{R}_t\}$ admits a rainbow matching?
\end{question}

%\section{The Proof of Theorem \ref{TH-AHC}}

\section{Rainbow matchings}

In this section we prove Theorem~\ref{TH-AHC}. First, we prove the
following lemma, which will serve as basis for our inductive proof of Theorem~\ref{TH-AHC}.

\begin{lemma}\label{Bipartite-rainbow}
Let $n>t>0$ be integers, and let $U_1, \ldots, U_k$ be pairwise disjoint sets with $|U_i|=n$ for $i\in [k]$.
 For each $i\in [t]$, let  $i_1,i_2\in [k]$ be distinct and let
 $F_i\subset U_{i_1}\times U_{i_2}$ such that $e(F_i)>(t-1)n$
for $i\in [t]$. Then $\{F_1,\ldots,F_t\}$ admits a rainbow matching.
\end{lemma}

\pf
Since $e(F_1)>(t-1)n$, there exists $x_1\in V(F_1)$ such that
$d_{F_1}(x_1)\geq t$.  Since $e(F_2)>(t-1)n$, there
exists $x_2\in V(F_2)$ such that $d_{F_2-x_1}(x_2)\geq t-1$.
 Suppose that we have chosen $x_{s-1}$, where
$2\le s-1\leq t-1$, such that $d_{F_{s-1}-\{x_1,\ldots,x_{s-2}\}}(x_{s-1})\geq t-(s-2)$.
Let $X_{s-1}=\{x_1,\ldots,x_{s-1}\}$,  $|X_{s-1}\cap U_{s_1}|=a$, and $|X_{s-1}\cap U_{s_2}|=b$. Then
$a+b\leq |X_{s-1}|=s-1$.

Without loss of generality, we may assume $a\geq b$. Then
\begin{align*}
e(F_s-X_{s-1})&>(t-1)n-(an+bn-ab)\\
&=\left(t-1-(a+b)\right)n+ab\\
&\geq (t-s)n+ab.
\end{align*}
Hence,   $F_{s}-X_{s-1}$ contains a vertex $x_s$ such
that
\begin{align*}
d_{F_S-X_{s-1}}(x_s)>\frac{(t-s)n+ab}{n-a} >t-s.
\end{align*}
%We choose from $V(F_{s+1})-X_s$ a vertex with degree more than $t-s-1$
%denoted by $x_{s+1}$.
Therefore, we obtain a sequence
$x_1,\ldots,x_t$ of distinct elements of $\bigcup_{i\in [k]}U_i$ such that
$d_{F_s-X_{s-1}}(x_{s})\ge t-(s-1)$ for $s\in [t]$, where $X_0=\emptyset$.

We now show that the desired rainbow matching exists by finding edges $e_s\in F_s$ in the order $s=t, \ldots, 1$.
Since $d_{F_t-X_{t-1}}(x_t)\ge 1$,  there exists $e_t\in F_t$ such that $e_t\cap X_{t-1}=\emptyset$
and $x_t\in e_t$. Suppose we have found pairwise disjoint edges $e_t, \ldots, e_{s+1}$
for some $s\in [t-1]$, such that, for $s+1\le j\le t$, $e_j\in F_j$, $e_j\cap X_{j-1}=\emptyset$, and $x_j\in e_j$.
Since $F_s$ is bipartite, $x_s$ is adjacent to
at most one vertex of each $e_j$, for $s+1\le j\le t$. Thus, since
$d_{F_{s}-X_{s-1}}(x_s)\ge t-(s-1)$, there exists  $e_s\in F_s$ such that $e_s\cap X_{s-1}=\emptyset$
and $x_s\in e_s$. Hence, by induction, there exist pairwise disjoint edges
$e_1,\ldots, e_t$ which form a rainbow matching for
$\{F_1,\ldots, F_t\}$. \qed

\bigskip

\noindent\textbf{Proof of Theorem \ref{TH-AHC}.}  %By Lemma \ref{shift}, we may assume that all $F_i's$ have been
%shifted with
%this enumeration.
We apply induction on $t+r$. Clearly, the assertion
holds for $t=1$. For $r=2$, the assertion follows from
Lemma~\ref{Bipartite-rainbow}. Therefore, we may assume $t\ge 2$ and
$r\ge 3$, and that the assertion holds with smaller $t+r$.

Suppose for $i\in [t]$, $|\{x\in V(\mathcal{F}_i): d_{{\cal
    F}_i}(x)>2(t-1)n^{r-2}\}|\ge t$.
Then there exist pairwise distinct $x_1,\ldots,x_t$ such that, for
$i\in [t]$, $x_i\in V(\mathcal{F}_i)$ and 
$d_{\mathcal{F}_i}(x_i)> 2(t-1)n^{r-2}$. Let $X:=\{x_1,\ldots,x_{t}\}$ and for $i\in [t]$,
let $\mathcal{F}_i':=\{S\ :\ S\subset V(\mathcal{F}_i)-(X-x_i)\mbox{ and } S\cup \{x_i\}\in {\cal F}_i\}$.
Then, for $i\in [t]$, $$|\mathcal{F}_i'| =d_{\mathcal{F}_i-(X-x_i)}(x_i)> 2(t-1)n^{r-2}-(t-1)n^{r-2}=(t-1)n^{r-2}.$$
%Let $\mathcal{F}_1':=N_{\mathcal{F}_i-(X-x_i)}(x_i)$ for $i\in [t]$.
By induction hypothesis, let  $\{e_1, \ldots, e_t\}$ be a rainbow matching
for $\{\mathcal{F}_1',\ldots,\mathcal{F}_t'\}$, with $e_i\in \mathcal{F}_i'$ for $i\in [t]$.
Clearly, $\{e_1\cup \{x_1\},\ldots,e_t\cup \{x_t\}\}$ is  a rainbow
matching for $\{\mathcal{F}_1,\ldots,\mathcal{F}_t\}$.

Hence, we may assume, without loss of generality, that $|\{x\in V(\mathcal{F}_t): d_{{\cal
    F}_t}(x)>2(t-1)n^{r-2}\}|\le t-1$.
By induction hypothesis, there exists a rainbow matching $M'$ for
$\{\mathcal{F}_1,\ldots,\mathcal{F}_{t-1}\}$.

Suppose $d_{\mathcal{F}_t}(x)\le (t-1)(r-1)n^{r-2}$ for all $x\in V(\mathcal{F}_t)$.
Then, since $r\ge 3$ and $|\{x\in V(\mathcal{F}_t): d_{{\cal
    F}_t}(x)>2(t-1)n^{r-2}\}|\le t-1$, the  number of edges in $\mathcal{F}_t$ intersecting $V(M')$ is
less than
\begin{align*}
&\quad (t-1)\left((t-1)(r-1)n^{r-2}\right)+\left((t-1)r-(t-1)\right)(2(t-1)n^{r-2})\\
%&=(t-1)^2(r-1)n^{r-2}+(t-1)(r-1)(2t-2)n^{r-2}\\
&= (t-1)^2(3r-3)n^{r-2}\\
&\le (t-1)n^{r-1}, \quad \mbox{since $n\geq(3k-3)(t-1)$ and $k\geq r\ge 3$}.
\end{align*}
So  there exists $e\in \mathcal{F}_t-V(M')$. Hence $M'\cup \{e\}$ is
rainbow matching for $\{{\cal F}_1,\ldots, {\cal F}_t\}$.

Therefore, we may assume that there exists  $x\in V(\mathcal{F}_t)$ such
that $d_{\mathcal{F}_t}(x)>(t-1)(r-1)n^{r-2}$.
%Let $E_i(x):=\{e\in \mathcal{F}_i\ :\ x\in e\}$ for $i\in [t]$; and
 For $i\in [t-1]$, let $\mathcal{F}_i''=\mathcal{F}_i- \{e\in \mathcal{F}_i\ :\ x\in e\}$. Then
\[
|\mathcal{F}_i''|> (t-1)n^{r-1}-n^{r-1}=(t-2)n^{r-1}.
\]
Hence, by induction bypothesis, there exists a rainbow matching $M$ for
$\{\mathcal{F}_1'',\ldots,\mathcal{F}_{t-1}''\}$.
Since $d_{\mathcal{F}_i''}(x)=0$ for $i\in [t-1]$, $x\notin V(M)$.

Since the number of edges  in $\mathcal{F}_t$  containing $x$ and intersecting $V(M)$ is at most
$(t-1)(r-1)n^{r-2}<d_{\mathcal{F}_t}(x)$,
there exists  $e\in \mathcal{F}_t-V(M)$. So $M\cup \{e\}$ gives the
desired rainbow matching for $\{F_1,\ldots,F_t\}$.\qed

\medskip

We now prove Conjecture~\ref{AH-conj} for the case when $t=n$.

\begin{proposition}\label{rainbow-nPM}
Let $k,n,t$ be  positive integers with $t\leq n$, and for $i\in [k]$,
let $W_i=\{jk+i\ :\ j\in [n-1]\cup \{0\}\}$.
For $i\in [n]$, let $\mathcal {F}_i\subset W_1\times\cdots \times W_k$ such that
$|\mathcal {F}_i|>(t-1)n^{k-1}$. Then there exist pairwise distinct
$i_1,\ldots, i_t$ from $[k]$
such that $\{\mathcal{F}_{i_1},\ldots,\mathcal{F}_{i_t}\}$ admits  a rainbow matching.
\end{proposition}

\pf Consider a permutation $\pi$ of $[kn]$, taken uniformly at random
from all permutations $\pi$ of $[kn]$ with the property that $\pi (W_i)=W_i$ for
all $i\in [k]$. For $i\in [n]$, let $X_i=1$ if
$\{\pi((i-1)k+1),\pi((i-1)k+2),\ldots,\pi(ik)\}\in \mathcal{F}_i$, and
let $X_i=0$ otherwise. Then
\begin{align*}
\mathbb{P}(X_i=1)=\frac{|\mathcal {F}_i|}{n^k}>\frac{(t-1)}{n}.
\end{align*}
Hence
\begin{align*}
\mathbb{E}\left(\sum_{i=1}^n X_i\right)=n
\mathbb{E}(X_i=1)>t-1.
\end{align*}
Therefore, there exist pairwise distinct $i_1,\ldots, i_t$ from $[k]$
such that $X_{i_j}=1$ for $j\in [t]$. Hence, $\{\mathcal{F}_{i_1},\ldots,\mathcal{F}_{i_t}\}$ admits  a rainbow matching.\qed

\medskip

Setting $t=n$ in Proposition \ref{rainbow-nPM}, we obtain the following
result on perfect matchings.
\begin{coro}
Let $n,k$ be  positive integers, and let
$\mathcal {F}_i\subset [n]^k$ for $i\in [n]$.  If $|\mathcal {F}_i|>(n-1)n^{k-1}$
for $i\in [n]$, then $\{\mathcal{F}_1,\ldots,\mathcal{F}_n\}$ admits a
rainbow matching.
\end{coro}

Setting $\mathcal{F}_1=\cdots=\mathcal{F}_n$ in Proposition \ref{rainbow-nPM}, we  obtain the following well-known result .
\begin{coro}\label{coro27}
Let $n,k$ be  positive integers, and let $H$ be a $k$-partite
$k$-graph with $V(H)=[n]^k$. If $e(H)>(t-1)n^{k-1}$, then $\nu(H)\ge t$.
\end{coro}

\section{Product type conditions}

In this section we prove Theorem~\ref{TH-HLS}.
First, we state a result of Matsumoto and Tokushige \cite{MN}.

\begin{lemma} \label{MTT}
Let $k_1,k_2,n$ be positive integers such that $n\geq \max\{2k_1,2k_2\}$,
and let $H_i\subset {[n]\choose k_i}$, $i\in [2]$, such that $e(H_1)e(H_2)>
{n-1\choose k_1-1}{n-1\choose k_2-1}$. Then  $\{H_1,H_2\}$ admits a rainbow
matching.
\end{lemma}

We use Lemma~\ref{MTT} as induction basis to prove the next result.

\begin{lemma}\label{lem-diff1}
Let $k,t,n$ be integers such that $t\ge 2$, $k_1\ge k_2\ge \ldots \ge k_t\ge 2$, and  $n\geq 9k_1^5t/k_2$.
%Let $\pi$ be a permutation of $[t]$, and let $\mathcal{F}_{i}\subseteq {[n] \choose k_{\pi(i)}}$, $i\in [t]$, such that
Let $\mathcal{F}_{i}\subset {[n] \choose k_{i}}$ for $i\in [t]$, such that
\begin{align*}%\label{F1F2}
 |\mathcal{F}_1|\ |\mathcal{F}_2|>{n-1\choose k_1-1}{n-1\choose k_2-1} \quad \mbox{when $$t=2},
\end{align*}
 and
\begin{align*}%\label{dffer-k_ii}
\prod_{i=1}^t |\mathcal{F}_i|>{n-1\choose k_1-1}{n-1\choose
  k_2-1}\prod_{i=3}^{t} {n\choose k_i} \quad \mbox{when $t\geq 3$.}
\end{align*}
Then $\{\mathcal{F}_1,\ldots,\mathcal{F}_t\}$ admits a rainbow matching.
 \end{lemma}

\pf  If $t=2$
then the assertion follows from Lemma~\ref{MTT}.  Thus, we may  assume
that $t\geq 3$ and the assertion holds with fewer than  $t$ families.
Let $s\in [t]$ such that
$$\frac{|\mathcal{F}_{s}|}{{n\choose k_{s}}}=\max\left\{\frac{|\mathcal{F}_{i}|}{{n\choose k_{i}}}:i\in [t]\right\}.$$

Since $|\mathcal{F}_{s}|\leq {n\choose k_{s}}$,
if $s\notin [2]$ then
\begin{align*}
\prod_{i\in [t]-\{s\}}|\mathcal{F}_{i}|&>  {n-1\choose k_1-1}{n-1\choose k_2-1}\prod_{i\in [t]-\{1,2,s\}} {n\choose k_{\pi(i)}},
\end{align*}
and if $s\in [2]$ then
\begin{align*}
\prod_{i\in [t]-\{s\}}|\mathcal{F}_i|&>  {n-1\choose k_1-1}{n-1\choose k_2-1}\prod_{i=3}^{t} {n\choose k_{i}}/{n\choose k_s}\\
&= \frac{k_{s}}{n}{n-1\choose k_{3-s}-1}\prod_{i=3}^{t} {n\choose k_{i}}\\
&=\frac{k_{s}}{k_{3}}{n-1\choose k_{3-s}-1}{n-1\choose
  k_{3}-1}\prod_{i=4}^{t} {n\choose k_{i}}\\
&\geq {n-1\choose k_{3-s}-1}{n-1\choose k_{3}-1}\prod_{i=4}^{t} {n\choose k_{i}} \quad \mbox{(since $k_1\geq k_2\geq k_{3}$)}.
\end{align*}

By induction hypothesis,
$\{\mathcal{F}_{1},\ldots,\mathcal{F}_{t}\}-\{{\cal F}_s\}$ admits a
rainbow matching, say $M$. Note that the number of edges  in
$\mathcal{F}_{s}$  intersecting
$V(M)$ is at most
\[
\left(\sum_{i\in [t]-\{s\}}k_{i}\right) {n-1\choose k_{s}-1}.
\]

Hence, if $|\mathcal{F}_s|>(\sum_{i\in [t]-\{s\}}k_{i}){n-1\choose k_{s}-1}$ then 
 there exists $e\in \mathcal{F}_{s}$ disjoint from $V(M)$. Thus $M\cup
\{e\}$ is  the desired rainbow matching for $\{\mathcal{F}_{1},\ldots,\mathcal{F}_{t}\}$.

So we may assume that $|\mathcal{F}_s|\le (\sum_{i\in
  [t]-\{s\}}k_{i}){n-1\choose k_{s}-1}$. Then, since $k_i\le k_1$ for
$i\in [t]$,
\[
|\mathcal{F}_s|\le k_1(t-1)\frac{k_{s}}{n}{n\choose k_{s}}\leq\frac{k_1^2(t-1)}{n}{n\choose k_{s}}.
\]
Therefore,
\begin{align*}%\label{frac:r1}
\prod_{i=1}^t \frac{|\mathcal{F}_{i})}{{n\choose k_{i}}}\le \left(\frac{|\mathcal{F}_{s}|}{{n\choose k_{s}}} \right)^t\leq \left(\frac{k_1^2(t-1)}{n}\right)^t<\left(\frac{k_1^2t}{n}\right)^t.
\end{align*}
On the other hand, by assumption of this lemma,
\begin{align*}%\label{frac:r2}
\prod_{i=1}^t \frac{|\mathcal{F}_{i}|}{{n\choose k_{i}}}>
  \frac{{n-1\choose k_1-1}{n-1\choose k_2-1}} {{n\choose k_1}{n\choose k_2}} = \frac{k_1k_2}{n^2}.
\end{align*}
Hence,
\begin{align*}%\label{k_t-k_1k_2}
\left(\frac{k_1^2t}{n}\right)^t>\frac{k_1k_2}{n^2}.
\end{align*}
Thus,
%By the logarithm of both ends of Inequality (\ref{k_t-k_1k_2}), we get
\begin{align*}\label{f(t)<0}
f(t):=t(\ln k_1^2+\ln t-\ln n)-\ln (k_1k_2)+2\ln n>0.
\end{align*}
However, the derivative $f'(t)=\ln( k_1^2)+\ln t-\ln n+1<0$, since
$n\geq 9k_1^5 t/k_2>3k_1^2t$. Thus $f(t)$ is a decreasing function.
Hence, since $3\leq t\leq k_2n/(9k_1^5)$ and $n\ge 9k_1^5t/k_2$, we have
\begin{align*}
f(t)\leq f(3)&=3(\ln k_1^2+\ln 3-\ln n)-\ln (k_1k_2)+2\ln n<0,
\end{align*}
a contradiction. \qed

\medskip

We need another lemma.

\begin{lemma}\label{lem-diff2}
Let $t$, $n$, $k_1, \ldots, k_t$ be positive integers and let
$\varepsilon>0$ be a small constant such that  $k_1\ge k_2\ge \ldots \ge k_t$, $\sum_{i=1}^t k_i\leq
n(1-\varepsilon)$, and $n$ is sufficiently large. For
$i\in [t]$, let $\mathcal{F}_i\subseteq {[n]\choose k_{i}}$ such that
\begin{align*}%\label{F1F2}
 |\mathcal{F}_1||\mathcal{F}_2|>{n-1\choose k_1-1}{n-1\choose k_2-1} \quad \mbox{when $$t=2},
\end{align*}
 and
\begin{align*}%\label{dffer-k_ii}
\prod_{i=1}^t |\mathcal{F}_i|>{n-1\choose k_1-1}{n-1\choose k_2-1}\prod_{i=3}^{t} {n\choose k_i} \quad \mbox{when $t\geq 3$.}
\end{align*}
%\begin{align}\label{dffer-k_i}
%\prod_{i=1}^t e(\mathcal{F}_i)>{n-1\choose k_1-1}{n-1\choose k_2-1}\prod_{i=3}^{t} {n\choose k_i}.
%\end{align}
Then $\{\mathcal{F}_1,\ldots,\mathcal{F}_t\}$ admits a rainbow matching.
 \end{lemma}

\pf  By Lemma \ref{lem-diff1}, we may assume that $t\geq
9k_1^5n/k_2$. Since $n$ is sufficiently large, we have $t\ge 3$. 
Let $s\in [t]$ such that $$\frac{|{\cal F}_s|}{{n\choose k_s}}
=\max \left\{\frac{|{\cal F}_i|}{{n\choose k_i}}: i\in [t]\right\}.$$
As induction hypothesis,  we may
assume that $M$ is a rainbow matching for $\{\mathcal{F}_1,\ldots,\mathcal{F}_{t}\}-\{{\cal F}_s\}$.
Then $|V(M)|= \sum_{i\in [t]-\{s\}}k_{i}.$
Again since $n$ is sufficiently large and $\sum_{i\in [t]}k_i\le n(1-\varepsilon)$, the number of edges  in
$\mathcal{F}_{s}$ intersecting $V(M)$ is at most
\begin{align*}
{n\choose k_{s}}-{n-\sum_{i\in [t]-\{s\}}k_{i}\choose k_{s}}\le {n\choose k_{s}}-{n\varepsilon\choose k_{s}}<{n\choose k_{s}}\left(1-\frac{1}{2}\varepsilon^{k_{s}}\right).
\end{align*}

If $|\mathcal{F}_s|>{n\choose k_{s}}-{n-\sum_{i\in [t]-\{s\}}k_{i}\choose k_{s}}$ then 
there exists  $e\in \mathcal{F}_{s}$ disjoint from
$V(M)$; so $M\cup \{e\}$ is the desired rainbow matching
for $\{\mathcal{F}_1,\ldots,\mathcal{F}_t\}$.

So assume  $|\mathcal{F}_s|\le {n\choose k_{s}}-{n-\sum_{i\in [t]-\{s\}}k_{i}\choose k_{s}}$.
By assumption of this lemma,
$$\prod_{i=1}^t \frac{|\mathcal{F}_i|} {{n\choose k_{i}}}>\frac{{n-1\choose k_1}{n-1\choose k_2}}{{n\choose k_1}{n\choose k_2}}=\frac{k_1k_2}{n^2}.$$
However, since $t\ge 9k_1^5n/k_2\ge 9k_1^4n$,
%\begin{align*}\label{cap-edge-up}
$$\frac{|\mathcal{F}_{s}|}{{n\choose k_{s}}}\geq \left(\prod_{i=1}^t \frac{|\mathcal{F}_i} {{n\choose k_{i}}}\right)^{1/t}>\left(\frac{k_1k_2}{n^2}\right)^{1/(9k_1^4n)}.$$
%\end{align*}
Since
\[
\lim_{n\rightarrow +\infty}\left(\frac{k_1k_2}{n^2}\right)^{1/(9k_1^4n)}=1
\]
and $n$ is sufficiently large,
\[
\frac{|\mathcal{F}_{s}|}{{n\choose k_{s}}}>1-\frac{1}{2}\varepsilon^{k_{s}}.
\]
This is a contradiction. \qed

\medskip

We also need the following inequality when $\sum_{i\in [t]}k_i\ge n(1-\varepsilon)$. .

\begin{lemma}\label{lem-inequ}
 Let $k,t,n$ be positive integers and $\varepsilon>0$ be a small constant, such that $n$ is sufficiently
 large, $k_1\geq k_2\geq \ldots\geq k_t\ge 2$, and
 \[
n(1-\varepsilon)\leq \sum_{i=1}^t k_i\leq n-n\left(\frac{8k_1\ln n}{n}\right)^{1/k_1}.\]
 Then
%\begin{align}\label{lem34-ieq}
\[
\left(\frac{k_1k_2}{n^2}\right)^{1/t}>1-\frac{1}{2}\left(1-\frac{\sum_{i=1}^t k_i}{n}\right)^{k_1}.
\]
%\end{align}
\end{lemma}

\pf Let  $m:=n-\sum_{i=1}^t k_{i}$. Then, by the assumption of this
lemma,  $$n \left(\frac{8k_1\ln n}{n}\right)^{1/k_1}\le   m\le n\varepsilon.$$
Moreover, $m=n-\sum_{i=1}^t k_{i}\ge n-tk_1$; so $n\varepsilon \ge n-tk_1$ and, hence,
$t\ge n(1-\varepsilon)/k_1.$
Since $n$ is large, we may assume $n^2>k_1k_2$. Hence
$$\frac{1}{t}\ln \left(\frac{k_1k_2}{n^2}\right)\ge \frac{k_1}{n(1-\varepsilon)}\ln \frac{k_1k_2}{n^2}>\frac{4k_1}{3n}\ln \frac{k_1k_2}{n^2},$$
where the last inequality holds as $\varepsilon$ is small (say
$\varepsilon <1/4$).

Note that the assertion of this lemma is equivalent to
\[
\left(\frac{k_1k_2}{n^2}\right)^{1/t}>1-\frac{1}{2}\left(\frac{m}{n}\right)^{k_1},
\]
which holds iff
$$\frac{1}{t}\ln\left(\frac{k_1k_2}{n^2}\right)>\ln \left( 1-\frac{1}{2}\left(\frac{m}{n}\right)^{k_1}\right).$$

Thus, since $n$ is large, it suffices to show
$$\frac{4k_1}{3n}\ln \frac{k_1k_2}{n^2}>-\frac{1}{3}\left(\frac{m}{n}\right)^{k_1}.$$
However, this follows from a straightforward calculation, using
$m\ge  n\left(\frac{8k_1\ln n}{n}\right)^{1/k_1}.$ \qed

\bigskip

\textbf{Proof of Theorem \ref{TH-HLS}.}
By Lemma~\ref{MTT}, we may assume $t\ge 3$. By Lemma \ref{lem-diff2}
and by assumption, we may assume that
 \[
n(1-\varepsilon)\leq \sum_{i=1}^t k_i\leq n-n\left(\frac{8k_1\ln n}{n}\right)^{1/k_1}.
\]
Let $s\in [t]$ such that $$\frac{|\mathcal{F}_{s}|}{{n\choose
    k_{s}}}=\max \left\{\frac{|\mathcal{F}_{i}|}{{n\choose
    k_{i}}}: i\in [t]\right\}.$$

As induction hypothesis, assume that $M$ is a rainbow matching for
$\{\mathcal{F}_1,\ldots,\mathcal{F}_{t}\}-\{{\cal F}_s\}$. Then
$|V(M)|=\sum_{i\in [t]-\{s\}} k_{i}$ and
the number of edges in $\mathcal{F}_s$ intersecting  $V(M)$ is at most
\begin{align*}
{n\choose k_s} -{n-\sum_{i\in [t]-\{s\}}k_{i}\choose k_{s}}.
\end{align*}

Note that  $$\frac{|\mathcal{F}_s|}{{n\choose k_{s}}}\geq \left(\prod_{i=1}^t \frac{|\mathcal{F}_i|}{{n\choose k_{i}}}\right)^{1/t}>
\frac{{n-1\choose k_1-1}{n-1\choose k_2-1}}{{n\choose k_1}{n\choose k_2}}=\left(\frac{k_1k_2}{n^2}\right)^{1/t}.$$
Hence, by Lemma \ref{lem-inequ} (and since $n$ is large),
$$\frac{|\mathcal{F}_s|}{{n\choose
    k_{s}}}>1-\frac{1}{2}\left(1-\frac{\sum_{i=1}^t
    k_i}{n}\right)^{k_{1}}
>1 -\frac{{n-\sum_{i=1}^t k_{i}\choose k_{s}}}{{n\choose k_{s}}}.$$
Therefore,
\[
|\mathcal{F}_{s}|>{n\choose k_{s}}-{n-\sum_{i=1}^t k_{i}\choose k_{s}}.
\]
So there exists $e\in \mathcal{F}_{s}$ such that $e\cap
V(M)=\emptyset$. Now $M\cup \{e\}$ is a rainbow matching for
$\{\mathcal{F}_1, \dots, \mathcal{F}_t\}$.
\qed

\section{Proof of Theorem \ref{HLS-EX12}}

Suppose the assertion of Theorem \ref{HLS-EX12} is false. We choose a counterexample so that $t$ is minimum and, subject to this,
$\sum_{i=1}^ t k_i$ is minimum. Clearly,  $t\ge 2$.

We  claim that  $k_i\ge 2$ for $i\in [t]$ and, for $i\in [t]$ and $v\in [n]$,
$d_{\mathcal{F}_i}(v)\le k(t-1){n-2\choose k_i-2}$.
For, suppose there exist $i\in [t]$ such that $k_i=1$ or 
$d_{\mathcal{F}_i}(v)>k(t-1){n-2\choose k_i-2}$ for some $v\in \mathcal{F}_i$. Note that, for $j\in [t]-\{i\}$,
\[
|\mathcal{F}_j-v|>{n\choose k_j}-{n-t+1\choose k_j}-{n-1\choose
  k_j-1}= {n-1\choose k_j}-{n-t+1\choose k_j}.
\]
Hence, $\{{\cal F}_1-v, \ldots, \mathcal{F}_t-v\}-\{{\cal F}_i-v\}$
admits a rainbow matching, say $M$. The number of edges in $\mathcal{F}_i$ containing $v$ and intersecting $V(M)$  is at most
\[
k(t-1){n-2\choose k_i-2}.
\]
So there exists $e\in {\cal F}_i$ such that $v\in e$ and $e\cap V(M)=\emptyset$. Therefore, $M\cup\{e\}$ is a rainbow matching for $\{\mathcal{F}_1,\ldots,\mathcal{F}_t\}$,
a contradiction.

Suppose for each $i\in [t]$,  $$\big|\left\{v\in [n]\ :\ d_{{\cal
    F}_i}(v)>2(t-1){n-2\choose k_i-2}\right\}\big|\ge t.$$
 Let $v_1,\ldots, v_t\in [n]$ be pairwise distinct  such that $d_{{\cal
    F}_i}(v_i)>2(t-1){n-2\choose k_i-2}$ for $i\in [t]$. Let $S_i=\{v_j:
j\in [t]-\{i\}\}$, $i\in [t]$. Then
\[
d_{\mathcal{F}_i-S_i}(v_i)>2(t-1){n-2\choose k_i-2} -(t-1){n-2\choose k_i-2} =(t-1){n-2\choose k_i-2}.
\]
 For $i\in [t]$, let
$$\mathcal{F}_i'=\left\{S: S\in {[n]-\{v_1,\ldots,v_t\}\choose k_i-1} \mbox{ and } S\cup \{v_i\}\in {\cal F}_i-S_i\right\}.$$
Then 
$$|{\cal F}_i'|=d_{\mathcal{F}_i-S_i}(v_i)>(t-1){n-2\choose k_i-2}.$$
So 
 $\{\mathcal{F}_1',\ldots,\mathcal{F}_t'\}$ admits a rainbow matching, say 
 $\{e_1,\ldots, e_t\}$ with $e_i\in \mathcal{F}_i'$ for $i\in [t]$.
Now $\{e_1\cup \{v_{1}\},\ldots,e_t\cup \{v_{t}\}\}$ is a
rainbow matching for $\{\mathcal{F}_1,\ldots,\mathcal{F}_t\}$, a contradiction.

Thus, without loss of generality, we may assume that
$$\big|\left\{v\in [n]\ :\ d_{{\cal
    F}_t}(v)>2(t-1){n-2\choose k_i-2}\right\}\big|<t.$$
 Let $M$ be a rainbow matching for $\{\mathcal{F}_1,\ldots,\mathcal{F}_{t-1}\}$.
Since  $\{\mathcal{F}_1,\ldots,\mathcal{F}_{t}\}$ admits no rainbow matching, every edge of $\mathcal{F}_t$ must intersect $V(M)$. Hence,
$|\mathcal{F}_t|$ is at most
\[
(k(t-1)-(t-1))\left(2(t-1){n-2\choose k_t-2}\right)+(t-1)\left(k(t-1){n-2\choose k_t-2}\right)=(3k-2)(t-1)^2{n-2\choose k_t-2}.
\]

On the other hand,
\begin{align*}
|\mathcal{F}_t|
&>{n\choose k_t}-{n-t+1\choose k_t}\\
&={n\choose k_t}-{n\choose k_t}\frac{(n-t+1)\cdots(n-t-k_t+2)}{n(n-1)\cdots(n-k_t+1)}\\
&>{n\choose k_t}\left(1-\left(1-\frac{t-1}{n}\right)^{k_t}\right)\\
&>{n\choose k_t}\left(\frac{k_t(t-1)}{n}-\frac{k_t^2(t-1)^2}{2n^2}\right)\\
&>{n\choose k_t}\frac{k_t(t-1)}{n}\left(1-\frac{1}{6k}\right) \quad
  (\mbox{since $n\geq 3k^2t$ and $k=\max\{k_i:i\in [t]\}$})\\
&={n-2\choose k_t-2}\frac{(n-1)(t-1)}{k_t-1}\left(1-\frac{1}{6k}\right)\\
&>{n-2\choose k_t-2}\frac{n(t-1)}{k_t}\left(1-\frac{1}{6k}\right),
\end{align*}

Therefore, $$(3k-2)(t-1)^2{n-2\choose k_t-2}>{n-2\choose k_t-2}\frac{n(t-1)}{k_t}\left(1-\frac{1}{6k}\right),$$
which implies
\[
n<k_t(3k-2)(t-1)\frac{6k}{6k-1}<3k^2t,
\]
a contradiction. \qed

\end{document}